\documentclass{ifacconf}

\usepackage{amsmath,amsfonts,amssymb}
\usepackage{graphicx,wrapfig}      
\usepackage{natbib}        		   

\newcommand{\be}{\begin{equation}}
\newcommand{\ee}{\end{equation}}
\newcommand{\bd}{\begin{displaymath}}
\newcommand{\ed}{\end{displaymath}}

\newcommand{\R}{\mathbb{R}}

\newcommand{\Htwo}{\mathcal{H}_{2}}

\newcommand{\w}{\omega}

\newcommand{\beqa}{\begin{eqnarray}}
\newcommand{\eeqa}{\end{eqnarray}}
\newcommand{\bena}{\begin{eqnarray*}}
\newcommand{\eena}{\end{eqnarray*}}
\newcommand{\beqn}{\begin{equation} \label}
\newcommand{\eeqn}{\end{equation}}
\newcommand{\bap}{\left[ \begin{array}}
\newcommand{\eap}{\end{array} \right]}
\newcommand{\mato}{\left( \begin{array}}
\newcommand{\matc}{\end{array} \right)}

\newcommand{\bea}{\begin{array}}
\newcommand{\eea}{\end{array}}

\newcommand{\smallmato}{\left( \begin{smallmatrix}}
\newcommand{\smallmatc}{\end{smallmatrix} \right)}

\begin{document}
\begin{frontmatter}

\title{SSD -- Software for Systems with Delays: \large{Reproducible Examples and Benchmarks on \\ Model Reduction and $\mathcal{H}_2$ Norm Computation}} 

\author[First]{Suat Gumussoy}
\address[First]{Siemens Technology, Princeton, NJ 08540 USA \\(e-mail: suat.gumussoy@siemens.com).}

\begin{abstract}
We present SSD -- \emph{Software for Systems with Delays}, a de novo MATLAB package for the analysis and model reduction of retarded time delay systems (RTDS). Underneath, our delay system object bridges RTDS representation and Linear Fractional Transformation (LFT) representation of MATLAB. This allows seamless use of many available visualizations of MATLAB. In addition, we implemented a set of key functionalities such as $\Htwo$ norm and system gramian computations, balanced realization and reduction by direct integral definitions and utilizing sparse computation. As a theoretical contribution, we extend the \emph{frequency-limited balanced reduction} to delay systems first time, propose a computational algorithm and give its implementation. We collected two set of benchmark problems on $\Htwo$ norm computation and model reduction. SSD is publicly available in GitHub\footnote{\texttt{https://github.com/gumussoysuat/ssd}}. Our reproducible paper and two benchmarks collections are shared as executable notebooks.
\end{abstract}


\end{frontmatter}

\section{Introduction}
Software packages are essential and practical tools for analysis and design of control systems. Considerable research effort is devoted to extend classical and modern control techniques to accommodate delays. Without being exhaustive and focusing only publicly available ones, we can group software packages as:
\begin{itemize}
\item Solution of delay differential equations, \cite{Hairer:95,Enright:97,Shampine:01}.
\item Stability and bifurcation analyses, \cite{Engelborghs:02,Breda:09,Wu:12,Avanessoff:13,Vhylidal:14,Breda:16}.
\item System norms and other metrics, \cite{Michiels:10,Gumussoy:10,Appeltans:19}.
\item Controller design: 
\begin{itemize}
\item Stabilization, \cite{Michiels:11b}.
\item Pole placement, \cite{Boussaada:21}.
\item PID design, \cite{Gumussoy:12,Appeltans:22}.
\item $\mathcal{H}_\infty$ control, \cite{Gumussoy:11}.
\end{itemize}
\end{itemize} 

This paper gives a guided tour to a new MATLAB package, Software for Delay Systems (SSD) focusing on model reduction and $\Htwo$ norm computation for retarded time delay systems (RTDSs). Our main contributions are
\begin{itemize}
\item allowing easy-to-access MATLAB's \emph{time and frequency domain visualizations} by bridging RTDS and MATLAB's LFT representation,
\item publicly available implementation of \emph{model reduction} and \emph{$\Htwo$ norm computation} using their direct definitions of the integral form and utilizing \emph{sparse computation},
\item as a \emph{theoretical contribution}, extension of the \emph{frequency-limited balanced reduction} first time for RTDS and a computational algorithm via integral expressions for delay systems,
\item collecting two sets of \emph{benchmark problems} on model reduction and $\Htwo$ norm computation and sharing them as two executable notebooks,
\item making our paper completely executable and reproducible to facilitate \emph{reproducible research}.
\end{itemize}

Our goals are two folds: First, to introduce publicly available functionalities as a baseline approach for comparison with the advanced techniques in model reduction and $\Htwo$ norm computation. Second, to facilitate the analysis of delay systems having small to mid-size state dimensions by SSD's easy-to-use interface.

The paper is organized as follows. First we define the delay system representation as an RTDS and show its use in the next section. We illustrate how time and frequency domain visualizations are used in Section~\ref{sec:visual}. The model reduction functionalities and $\Htwo$ norm computation are overviewed in Section~\ref{sec:modelreduction}. We summarized the collected benchmark problems in Section~\ref{sec:benchmarks}. We give details on the computational aspects of our MATLAB package in Section~\ref{sec:computation}. We end our paper with concluding remarks and some future directions.

\section{Delay System Definition} \label{sec:howto}
SSD constructs a retarded time-delay system as
\beqa \label{eq:ssd}
\dot{x}(t) &=& \sum_{i=0}^{m_A} \ A_i\ x(t-h^A_i)+ \sum_{i=0}^{m_B} \ B_i\ u(t-h^B_i), \\
\nonumber y(t) &=& \sum_{i=0}^{m_C} \ C_i \ x(t-h^C_i)+ \sum_{i=0}^{m_D} \ D_i \ u(t-h^D_i).
\eeqa The system matrices are $A_i\in\R^{n\times n}$, $B_i\in\R^{n\times n_u}$, $C_i\in\R^{n_y\times n}$ and $D_i\in\R^{n_y \times n_u}$. The number of time delays are $m_A$, $m_B$ for states and inputs in the state equation and $m_C$, $m_D$ for states and inputs in the output equation. By convention, $h_0^A=h_0^B=h_0^C=h_0^D=0$. All delays are non-negative real numbers.

The ssd object is constructed using the function \verb"ssd" as
\begin{verbatim}
>> sys = ssd(A,hA,B,hB,C,hC,D,hD);
\end{verbatim} where $A$, $B$, $C$, and $D$ matrices are $3$-dimensional matrices and $hA$, $hB$, $hC$, and $hD$ are unique, increasing non-negative vectors as shown below in Figure~\ref{fig:ssd}.
\begin{figure}[h]
\centerline{\includegraphics[width=0.9\linewidth]{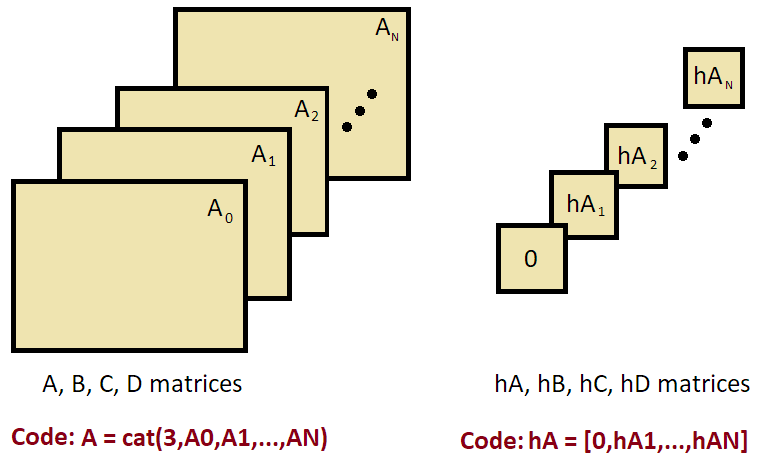}}
\caption{System matrices and delays in ssd object.\label{fig:ssd}}
\end{figure}

We borrow the following delay system from \cite{Jarlebring:13} as a motivational example,
\beqa
\nonumber \dot{x}(t)&=&\left[ {\begin{array}{cc}
    -2 & -1 \\
    -3/2 & -1/2 \\
  \end{array} } \right] x(t) + 
\left[ {\begin{array}{cc}
    0 & 1/2 \\
    1 & 0 \\
  \end{array} } \right] x(t-1) + 
\left[ {\begin{array}{c}
    1 \\
    -1 \\
  \end{array} } \right] u(t), \\
\nonumber y(t)&=&\left[ {\begin{array}{cc}
    2 & 0.2 \\
  \end{array} } \right] x(t).
\eeqa We first define the system matrices as follows, 
\begin{verbatim}
>> A0 = [-2 -1; -3/2 -1/2];
>> A1 = [0 1/2; 1 0];
>> A = cat(3,A0,A1);
>> B = [1; -1];
>> C = [2 0.2];
\end{verbatim} Then we create the ssd object by
\begin{verbatim}
>> sys = ssd(A,[0 1],B,0,C,0)

sys = 

  ssd with properties:

       A: [2×2×2 double]
      hA: [0 1]
       B: [2×1 double]
      hB: 0
       C: [2 0.2000]
      hC: 0
       D: 0
      hD: 0
    name: []
\end{verbatim}
Note that there is no $D$ matrix in the above example. Therefore, the matrix and its delay argument are not included. Alternatively, one can set them to empty vector, $\left[\ \right]$ or enter zero matrix with appropriate dimensions and zero delay.

The interface of ssd is a natural extension of state-space object, \verb"ss(A,B,C,D)", of MATLAB by introducing extra delays per system matrix. While the ssd object keeps the delay system structure, it also constructs the MATLAB's LFT-based ss object with delays behind the scenes to leverage pre-defined functionalities of MATLAB as we illustrate in the next section.

\section{Time / Frequency Domain Visualizations} \label{sec:visual}
Time and frequency domain plots gives additional insights for the analysis and design of delay systems. The ssd object interfaces with MATLAB plotting functionality and enables frequency domain plots, \verb"bode", \verb"bodemag", \verb"sigma", \verb"nyquist" and the time domain plot, \verb"step" as in MATLAB including multiple system support.

Continuing with the same system matrices of previous example, \verb"sys", we define two additional systems, one with a state delay of $0.5$ and the other one with no delay.
\begin{verbatim}
>> sys1 = ssd(A,[0 0.5],B,0,C,0);
>> sys2 = ssd(sum(A,3),0,B,0,C,0);
>> sys2.name='no delay';
>> step(sys,'r',sys1,'g.-',sys2,'b--');
\end{verbatim}
\begin{figure}[h]
\centering
\begin{minipage}[t][][b]{0.25\textwidth}
\centering
  \includegraphics[width=\linewidth]{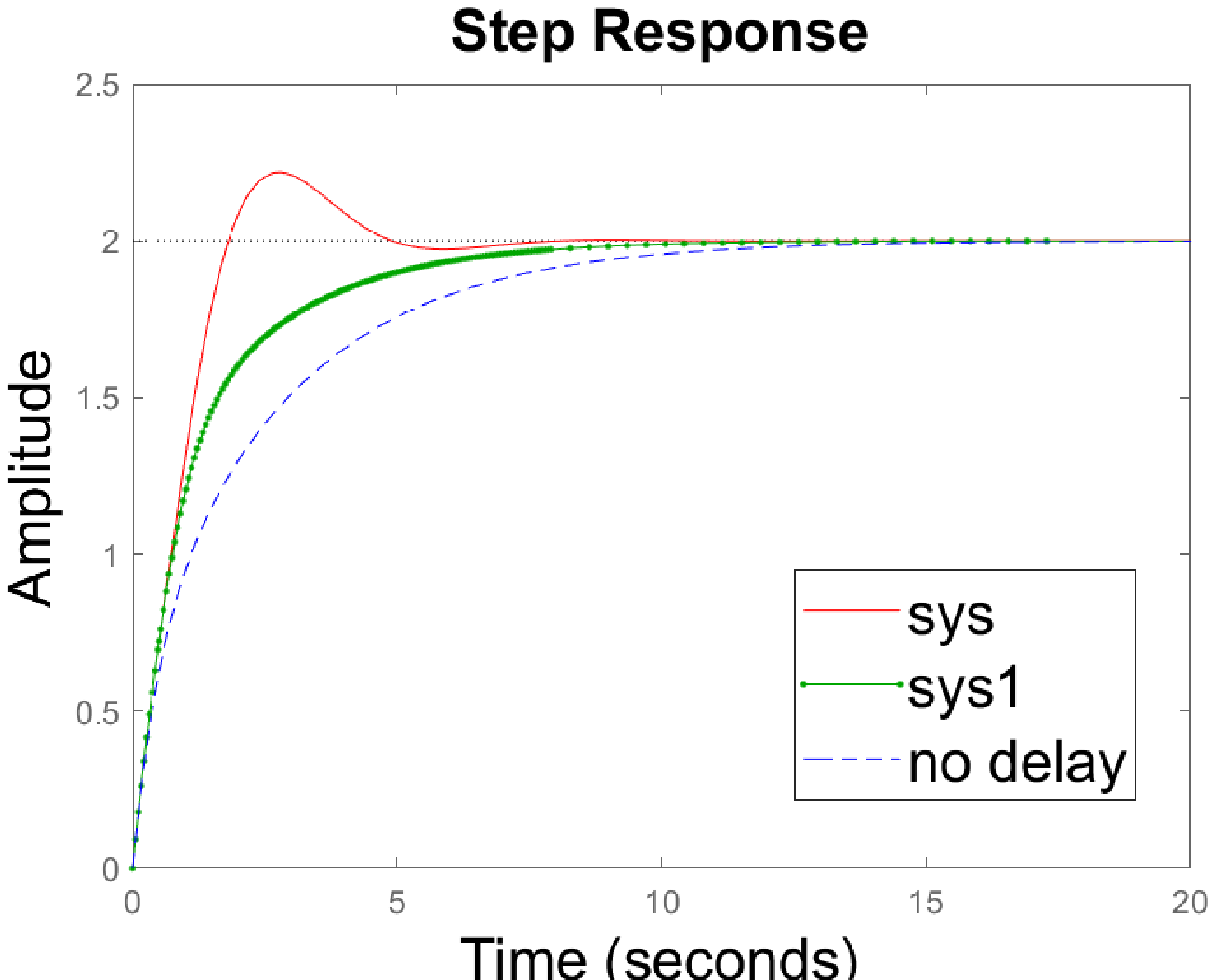}
\end{minipage}%
\begin{minipage}[t][][b]{0.25\textwidth}
\centering
  \includegraphics[width=\linewidth]{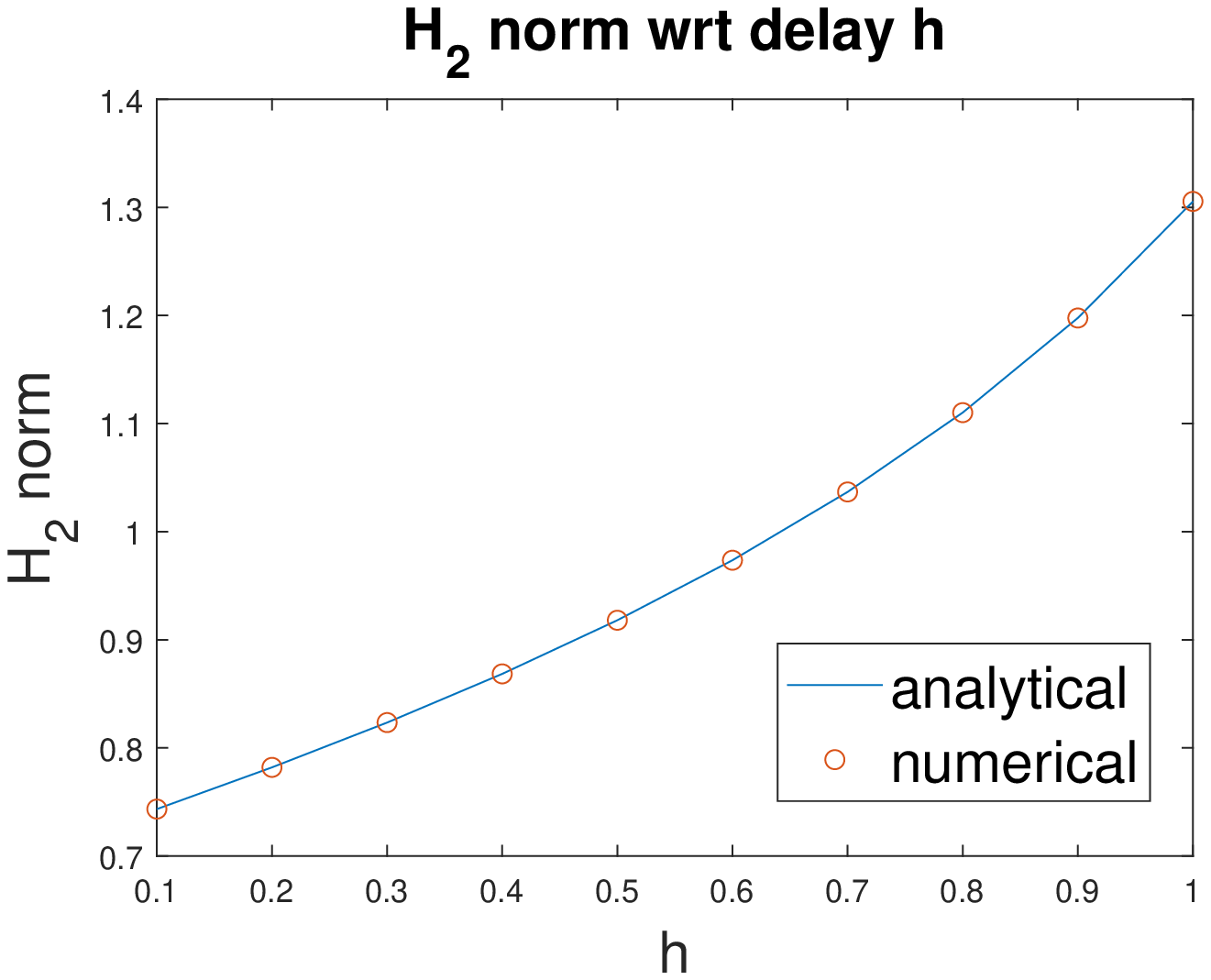}
\end{minipage} 
\par
\begin{minipage}[t]{.23\textwidth}
  \centering
  \caption{Step responses of multiple systems.\label{fig:step}}
\end{minipage}\hspace{3mm}
\begin{minipage}[t]{.23\textwidth}  
  \centering
  \caption{Analytical and numerical $\Htwo$ norm values. \label{fig:h2exact}}
\end{minipage}  
\end{figure}

Note that \verb"sys2" is a standard state-space system with no delay. As seen in Figure~\ref{fig:step}, SSD shows the system names and allows the user to modify the names as needed.

\section{$\Htwo$ Norm and Model Reduction} \label{sec:modelreduction}
SSD provides a set of functionality on $\Htwo$ norm computation and model reduction. In next two sections, we introduce the basic definitions, outline our direct computational approach and illustrate their use with the example codes.

\subsection{$\Htwo$ Norm Computation}
Assuming the exponential stability of the delay system, its $\Htwo$ norm is defined in frequency domain as,
\be \label{def:h2norm}
\|G\|_2: = \sqrt{\frac{1}{2\pi} \int_{-\infty}^\infty trace\left(G(j\omega)^*G(j\omega)\right)dw}
\ee where $G(s)=\mathcal{C}(s)\left(sI-\mathcal{A}(s)\right)^{-1}\mathcal{B}(s)$ with the auxiliary terms, $\mathcal{A}(s):=\sum_{i=0}^{m_A} A_i e^{-j\w h_i^A}$, $\mathcal{B}(s):=\sum_{i=0}^{m_B} B_i e^{-j\w h_i^B}$, and $\mathcal{C}(s):=\sum_{i=0}^{m_C} C_i e^{-j\w h_i^C}$. By definition, $D$ matrix is zero for a well-defined $\Htwo$ norm.

We borrow the hot shower problem from \cite{Jarlebring:11} (included in \verb"benchmarks_h2norm.mlx" at GitHub as the first benchmark). The system equations are
\[
\dot{x}=-ax(t-h)+bu(t),\quad y(t)=cx(t).
\] This example admits a closed form for the $\Htwo$ norm as, $\|G\|_2^2=\frac{c^2b^2}{2a}\frac{cos(ah)}{1-sin(ah)}$. Set the system parameters as $a=b=c=1$, then the following code computes the analytical and approximate $\Htwo$ norms.
\begin{verbatim}
>> a = 1; b = 1; c = 1;
>> h = 0.1:.1:1;
>> h2 = ((c*b)^2/2/a)*(cos(a*h)./(1-sin(a*h)));
>> for ct=1:length(h)
>>    sys = ssd(cat(3,0,-a),[0 h(ct)],1,0,1,0);
>>    happrox(ct) = h2norm(sys);
>> end
>> plot(h,sqrt(h2),'o',h,happrox);
\end{verbatim}
Figure~\ref{fig:h2exact} shows the comparison. The values are same within numerical tolerances.

SSD utilizes \emph{sparse computation} inside the integrals for $\Htwo$ norm computation (\ref{def:h2norm}). For delay systems with large dimensions, enabling \verb"sparse" option inside \verb"ssdoptions" can speed up the computation considerably. Computing $\Htwo$ norm without and with the sparse option enabled for the heated rod example \verb"sys4" of order $1000$ (included in our benchmark problems) results in $50$ times speed up.
\begin{verbatim}
>> tic; h2 = h2norm(sys4); time=toc

time = 400.7489
  
>> opt = ssdoptions('sparse',true);
>> tic; h2_sparse = h2norm(sys4,opt); time=toc

time = 8.3189
\end{verbatim}

\subsection{Balanced Reduction}
SSD implements one of the standard model reduction techniques, the balanced reduction (a survey of model reduction techniques is given in \cite{Antoulas:02}). In balanced reduction, a transformation matrix is computed where the transformed states are ordered from easy to difficult according to the measure of how difficult to observe and reach them simultaneously. Such transformation matrix is the {\it balancing transformation} and the transformed system is the {\it balanced realization}. The measure of reachability and observability for delay systems is defined in terms of {\it system gramians} in frequency domain as,
\beqa
W_c&=&\frac{1}{2\pi}\int_{-\infty}^\infty G_{\mathcal{A}\mathcal{B}}(j\w)G_{\mathcal{A}\mathcal{B}}^*(j\w) d\w, \label{def:Wc} \\
W_o&=&\frac{1}{2\pi}\int_{-\infty}^\infty G_{\mathcal{C}\mathcal{A}}^*(j\w)G_{\mathcal{C}\mathcal{A}}(j\w) d\w \label{def:Wo}.
\eeqa where $G_{\mathcal{A}\mathcal{B}}(s)=\left(s I-\mathcal{A}(s)\right)^{-1}\mathcal{B}(s)$ and $G_{\mathcal{C}\mathcal{A}}(s)=\mathcal{C}(s)\left(s I-\mathcal{A}(s)\right)^{-1}$.

Note that we define the gramians in the sense of position balancing, \cite{Jarlebring:13}, p.$160$. For convenience, we drop the term {\it position} in our paper.

The function \verb"gram" computes the system gramians. It calculates the system gramians by the direct definitions (\ref{def:Wc}) and (\ref{def:Wo}) and the \verb"integral" function of MATLAB. For our example system, \verb"sys", the system gramians are
\begin{verbatim}
>> [Wc,Wo] = gram(sys,'co')

Wc =

    0.9273   -1.7426
   -1.7426    3.6292

Wo =

    1.2674   -0.4129
   -0.4129    0.3674
\end{verbatim} which are equal to the values given in \cite{Jarlebring:13}, p.$160$. The sparse computation is also available for system gramians. It can be set from \verb"ssdoptions" and passed to the \verb"gram" function as an argument.

\begin{figure}[h]
\centering
\begin{minipage}[t][][b]{0.25\textwidth}
\centering
  \includegraphics[width=\linewidth]{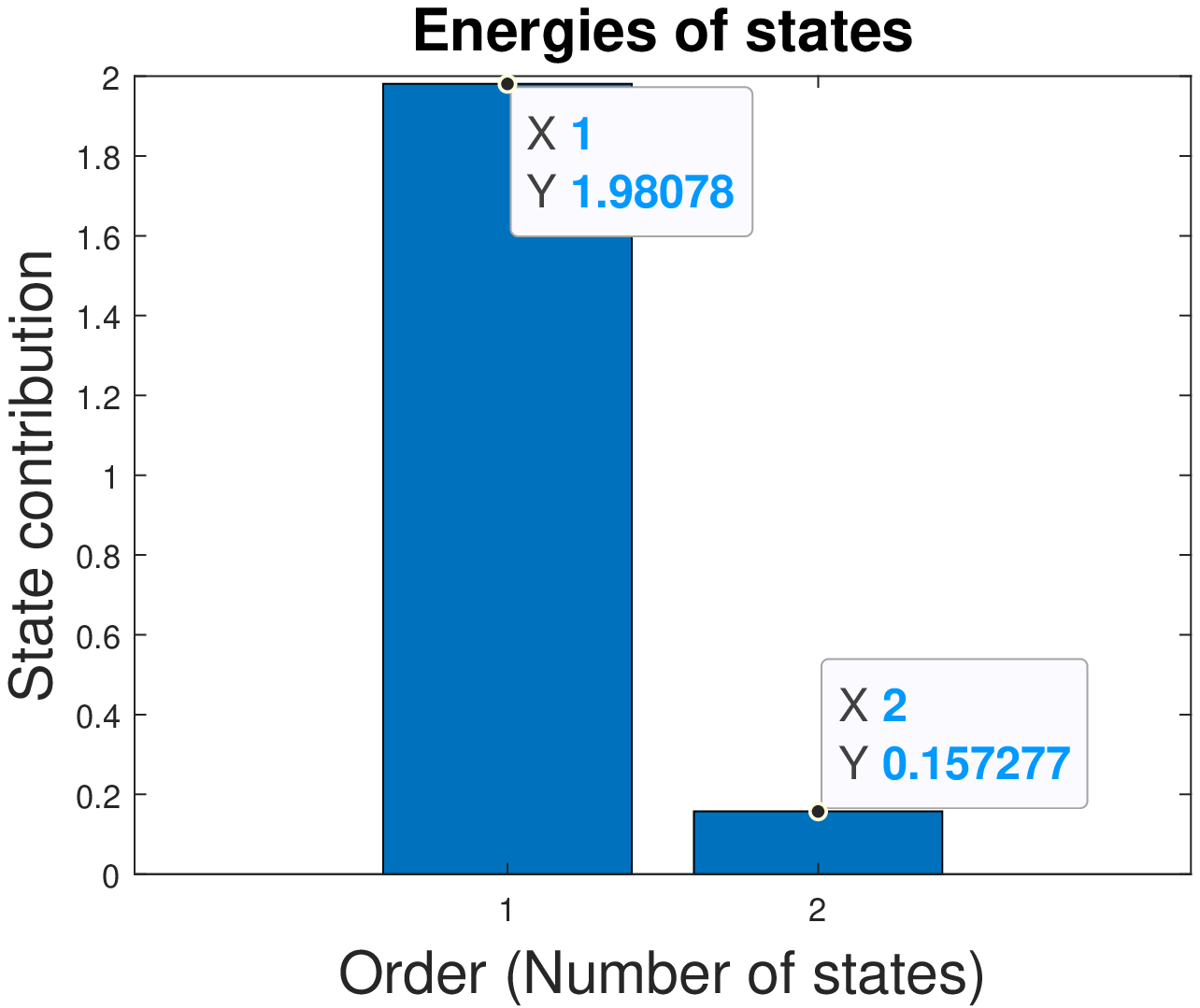}
\end{minipage}%
\begin{minipage}[t][][b]{0.25\textwidth}
\centering
  \includegraphics[width=\linewidth]{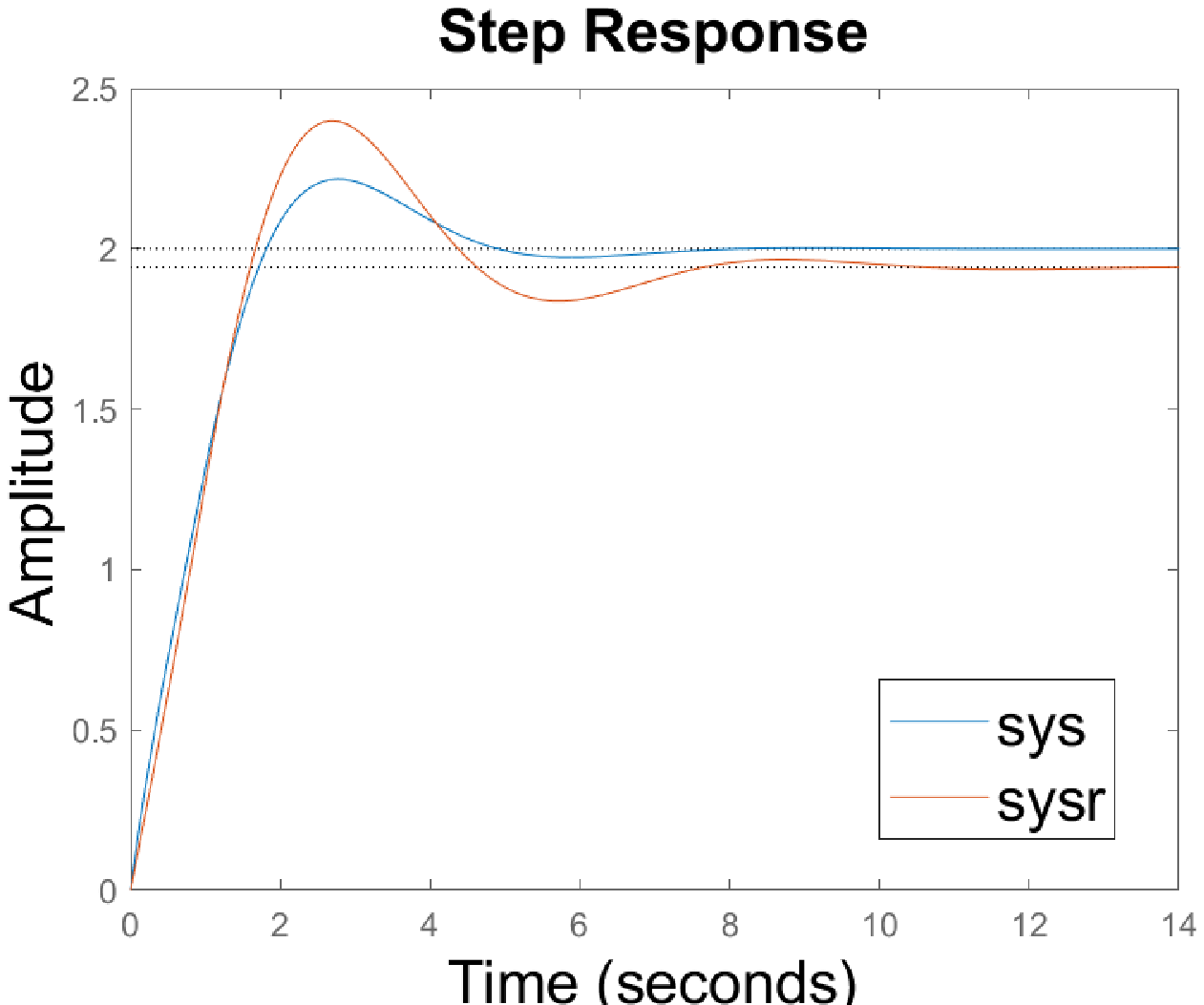}
\end{minipage} 
\par
\begin{minipage}[t]{.23\textwidth}
  \centering
  \caption{Energy contributions of the states.\label{fig:info_plot}}
\end{minipage}\hspace{3mm}
\begin{minipage}[t]{.23\textwidth}  
  \centering
  \caption{Step responses of \texttt{sys} and \texttt{sysr}.\label{fig:balred_step_plot}}
\end{minipage}  
\end{figure}

The function \verb"balreal" computes the balanced realization. It first computes the system gramians as above, then finds the balancing transformation as in Algorithm $1$ in \cite{Jarlebring:13}, p.$160$ without any truncation and applies to the delay system.
\begin{verbatim}
>> [sysb,info] = balreal(sys);
>> plot(info)
\end{verbatim} where \verb"sysb" is a balanced realization of \verb"sys" with the same system matrices given in \cite{Jarlebring:13}, p. $162$. One can see that the state contributions by plotting the \verb"info" object. Figure~\ref{fig:info_plot} shows that the first transformed state's contribution is major.

Based on the energy contributions of the states, the balanced reduction eliminates the states of the balanced realization upto a desired order defined by the user. The function \verb"balred" performs this reduction. Figure~\ref{fig:info_plot} suggests that the first-order model may approximate the \verb"sys" well. We reduce the system and compare the step responses of the original and the reduced systems by
\begin{verbatim}
>> sysr = balred(sys,1);
>> step(sys,sysr);
\end{verbatim}
Figure~\ref{fig:balred_step_plot} shows that the first-order model captures most of the dynamics as expected.

As a theoretical contribution, to the best of the author's knowledge, SSD extends the {\it frequency-limited balanced reduction} to delay systems first time where the system gramians are computed over the frequency intervals of interest (see \cite{Gawronski:90} for standard case). The resulting reduced system approximates the dynamics better over the desired frequency range with less number of states. 

Let \verb"sys3" represent the delay system Example $3$ in \cite{Jarlebring:11} (included as the third benchmark in \verb"benchmarks_h2norm.mlx"). The frequency interval is defined in \verb"FreqInt" field of  \verb"ssdoptions" object.
\begin{verbatim}
>> sys3r = balred(sys3,2);
>> opt = ssdoptions('FreqInt',[20 Inf]);
>> sys3rf = balred(sys3,2,opt);
>> sigma(sys3,sys3r,sys3rf);
\end{verbatim}

\begin{wrapfigure}{l}{0.25\textwidth}
\includegraphics[width=\linewidth]{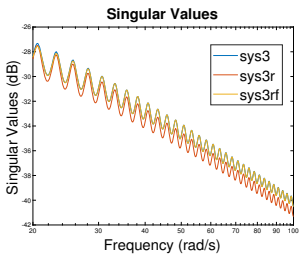}
\caption{Reduced systems \texttt{sys3r} and \texttt{sys3rf} over $[20, \infty)$.\label{fig:balred_freq}}
\end{wrapfigure} The zoomed version of the sigma plot, Figure~\ref{fig:balred_freq}, shows that the frequency-limited balanced reduction approximates better compared to the standard reduction over the focused frequency range.

As a final remark, the object \verb"info" stores the gramian related information and the full balancing transformation matrices. The function \verb"balred" computes the reduced orders very fast when \verb"info" is provided since it does not recompute gramian information, i.e.,
\begin{verbatim}
>> tic; [sys3r,info3] = balred(sys3,2); time=toc
>> tic; sys3r1 = balred(sys3,4,info3); time=toc

time = 3.1990

time = 0.0041
\end{verbatim}

\section{Benchmark Problems} \label{sec:benchmarks}
We surveyed the literature and collected two sets of benchmark problems for the model reduction and the $\Htwo$ norm computation. We prepared two executable notebooks in which we define the benchmark problems and show how to use SSD on them in GitHub repository.

\subsection{Model Reduction} \label{sec:bench_mr}
There are five benchmark problems provided in the executable notebook, \verb"benchmarks_model_reduction.mlx". We used the model reduction functionality of SSD on the selected benchmarks:
\begin{enumerate}
\item \textbf{The heated rod (HR$2$) Example in \cite{Michiels:11a}:} A $100^{\textrm{th}}$-order model with a single delay derived from the discretization of heat equations describing the temperature in a rod controlled with distributed delayed feedback.
\item \textbf{Mass-spring (MS) system in \cite{Saadvandi:12}:} A $1000^{\textrm{th}}$-order coupled mass-spring system with dampers and feedback controls with delays.
\item \textbf{Platoon of eight vehicles (P$8$V) in \cite{Scarciotti:14}:} A $23^{\textrm{rd}}$-order model describing the problem of controlling a group of vehicles tightly spaced following a leader, all moving in longitudinal direction.
\item \textbf{Example $1$ in \cite{Lordejani:20}:} A $6^{\textrm{th}}$-order synthetic model with one state and one output delay.
\item \textbf{Second-order system with proportional damping example (SOSPD) in \cite{Jiang:19}:} A $2000^{\textrm{th}}$-order model with a single state delay.
\end{enumerate}

Table~\ref{tbl:mrbenchmarks} shows the state, output, input dimensions and the number of delays for all system matrices. Since the maximum delay may affect the complexity of the method, this information is included in the table. 

The examples in \cite{Saadvandi:12} and \cite{Jiang:19} are second-order mechanical systems and the orders of their first-order equivalent models are twice of their dimensions. The benchmark problems are single-input-single-output (SISO) systems and mostly single state delay in the state equation is considered.

\begin{table}[h]
\begin{center}
\begin{tabular}{|c|c|c|c|}
  \hline
  \hline
  Ex. & Dimensions        & \# Delays & Max. \\
   & ($n$,$n_y$,$n_u$) & ($m_A$,$m_B$,$m_C$,$m_D$) & Delay\\
  \hline
  HR$2$ & $(100,1,1)$ & $(1,0,0,0)$ & $1$ \\  
  MS & $(1000^\dag,1,1)$ & $(1,0,0,0)$ & $2$ \\  
  P$8$V & $(23,1,1)$ & $(1,0,0,0)$ & $0.005$ \\    
  Ex.$1$ & $(6,1,1)$ & $(1,0,1,0)$ & $1.6$ \\
  SOSPD & $(2000^\dag,1,1)$ & $(1,0,0,0)$ & $1$ \\  
  \hline
  \hline
\end{tabular}
\newline
$\dag$ The dimension of the second-order system.
\caption{Model reduction benchmark problems} \label{tbl:mrbenchmarks}
\end{center}
\end{table}

\begin{figure}[h]
\centering
\begin{minipage}[t][][b]{0.25\textwidth}
\centering
  \includegraphics[width=\linewidth]{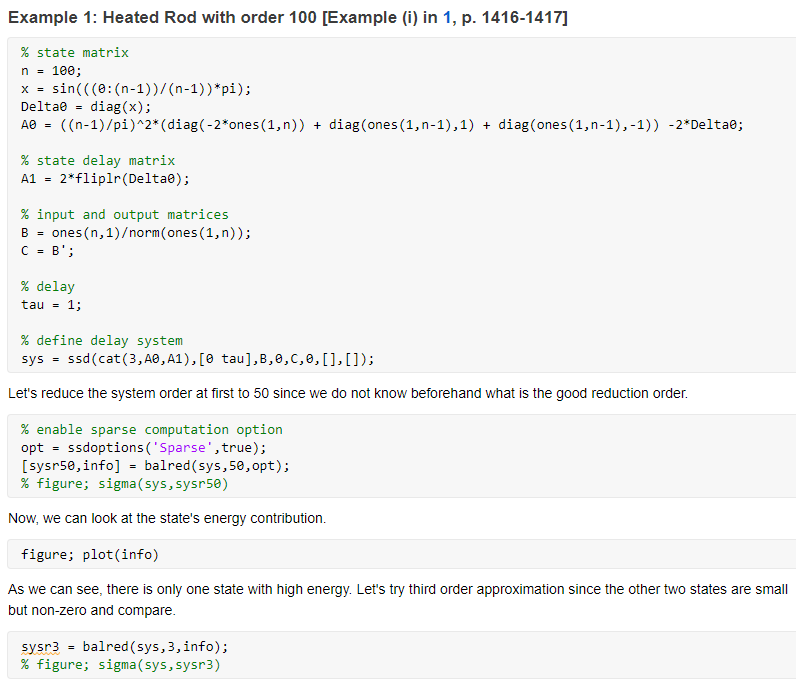}
\end{minipage}%
\begin{minipage}[t][][b]{0.25\textwidth}
\centering
  \includegraphics[width=\linewidth]{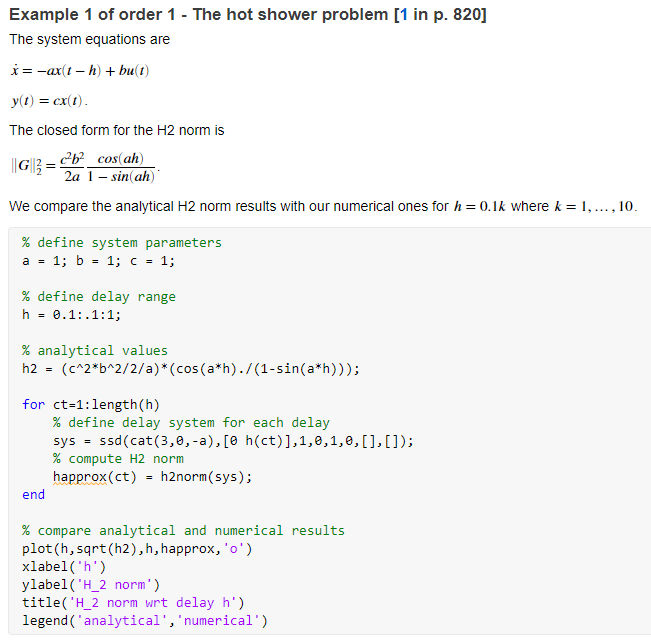}
\end{minipage} 
\par
\begin{minipage}[t]{.23\textwidth}
  \centering
  \caption{A snapshot of the heated rod example.\label{fig:mrbench}}
\end{minipage}\hspace{3mm}
\begin{minipage}[t]{.23\textwidth}  
  \centering
  \caption{A snapshot of the hot shower example.\label{fig:h2bench}}
\end{minipage}  
\end{figure}

A snapshot of the heated rod example from our notebook can be seen in Figure~\ref{fig:mrbench}.

\subsection{$\Htwo$ Norm Computation} \label{sec:bench_h2}
There are six benchmark problems for $\Htwo$ norm computation collected in our notebook, \verb"benchmarks_h2norm.mlx":
\begin{enumerate}
\item \textbf{The hot shower (HS) problem in \cite{Jarlebring:11}, Ex.$1$:} A first-order state-space model whose $\Htwo$ norm can be computed analytically. Its delay can be set by the user.
\item \textbf{Example $2.b$ in \cite{Jarlebring:11}, Ex.$2$b:} A $3^{\textrm{rd}}$-order synthetic model with two state delays.
\item \textbf{Example $3$ in \cite{Jarlebring:11}, Ex.$3$:} A $9^{\textrm{th}}$-order model with two state-delays obtained from discretization of partial differential equations with delays.
\item \textbf{The heat exchanger (HE) example in \cite{Michiels:19}:} A $5^{\textrm{th}}$-order model with seven delays for a heat exchanger for which the controller based on a combination of static state feedback and proportional integral control.
\item \textbf{The heated rod (HR$2$) Example in \cite{Peeters:13}:} A $100^{\textrm{th}}$-order model with a single delay, same as the previous example in the model reduction benchmarks.
\item \textbf{The heated rod (HR$4$) Example in \cite{Michiels:19}:} A $10000^{\textrm{th}}$-order model with a single delay, same as the previous example with larger size.
\end{enumerate}

\begin{table}[h]
\begin{center}
\begin{tabular}{|c|c|c|c|}
  \hline
  \hline
  Ref. & Dimensions        & \# Delays & Max. \\
   & ($n$,$n_y$,$n_u$) & ($m_A$,$m_B$,$m_C$,$m_D$) & Delay\\
  \hline
  HS & $(1,1,1)$ & $(1,0,0,0)$ & User-defined \\
  Ex.$2$b & $(3,1,1)$ & $(2,0,0,0)$ & $1$ \\  
  Ex.$3$ & $(9,1,1)$ & $(2,0,0,0)$ & $3$ \\    
  HE & $(5,1,5)$ & $(7,0,0,0)$ & $40$ \\
  HR$2$ & $(100,1,1)$ & $(1,0,0,0)$ & $1$ \\  
  HR$4$ & $(10000,1,1)$ & $(1,0,0,0)$ & $1$ \\
  \hline
  \hline
\end{tabular}
\newline
\caption{$\Htwo$ benchmark problems} \label{tbl:H2benchmarks}
\end{center}
\end{table}

Table~\ref{tbl:H2benchmarks} summarizes the dimensions and the number of delays for all system matrices. In general, the examples in the literature has one or two state delays and no delays at input and output matrices excluding $D$ matrices which must be zero for well-posedness of the norm. Most examples are SISO systems. Finally, we share a snapshot from the executable notebook for the hot shower example in Figure~\ref{fig:h2bench}.

\section{Computational Aspects} \label{sec:computation}
For experimental\footnote{All experiments are performed on HP Z6 G4 Workstation with Intel Xeon Silver 4114 2.2GHz CPU, 64GB RAM.} testing of the computational complexity of balanced reduction and $\Htwo$ norm computation, we randomly generate delay systems with the following properties:
\begin{itemize}
\item Input and output sizes are $3$, $n_y=n_u=3$,
\item There are $3$ delays, $m_A=m_B=m_C=3$,
\item All delays are uniformly selected from $(0, 1)$,
\item All matrices are dense whose elements uniformly selected from $[-3, 3]$ and in addition, $3I$ subtracted from $A_0$ matrix to make it stable.
\item $D$ matrices are set to zero since they are not needed for $\Htwo$ norm computation and balanced reduction.
\end{itemize}  

We validated that the generated delay systems have rich dynamics. We measured the computation time of \verb"h2norm" and \verb"gram" (computing both gramians) functions since the latter function is the main computational bottleneck of \verb"balreal" and \verb"balred". The number of states varies as $1-10$, $10-100$ by $10$ increments and $100-500$ for by $100$ increments (The upper bound is $400$ for model reduction). For each step, we created $5$ delay systems, computed the average of computational times and the standard deviation.

\begin{figure}[h]
\centering
\begin{minipage}[b]{0.25\textwidth}
\centering
  \includegraphics[width=\linewidth]{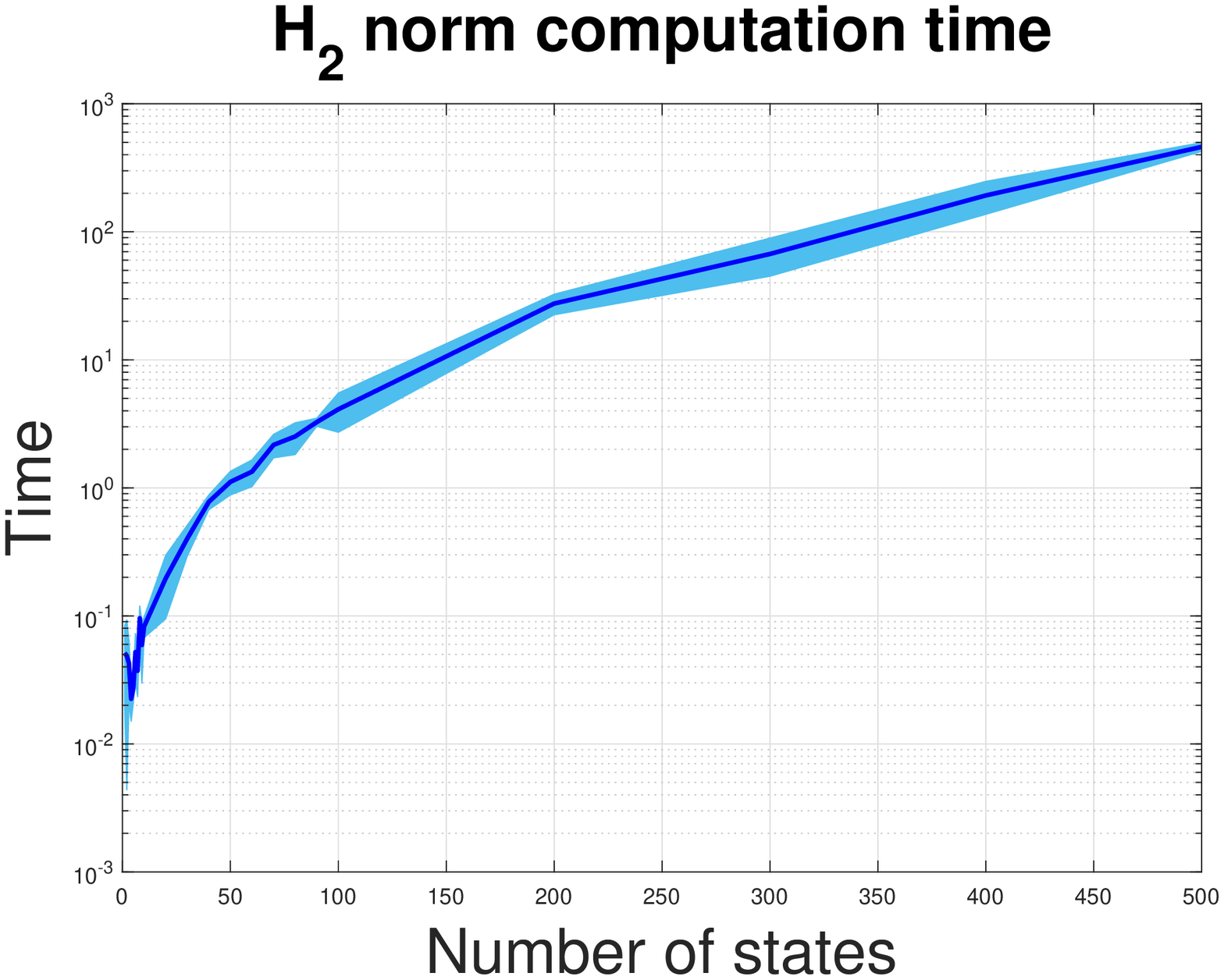}
\end{minipage}%
\begin{minipage}[b]{0.25\textwidth}
\centering
  \includegraphics[width=\linewidth]{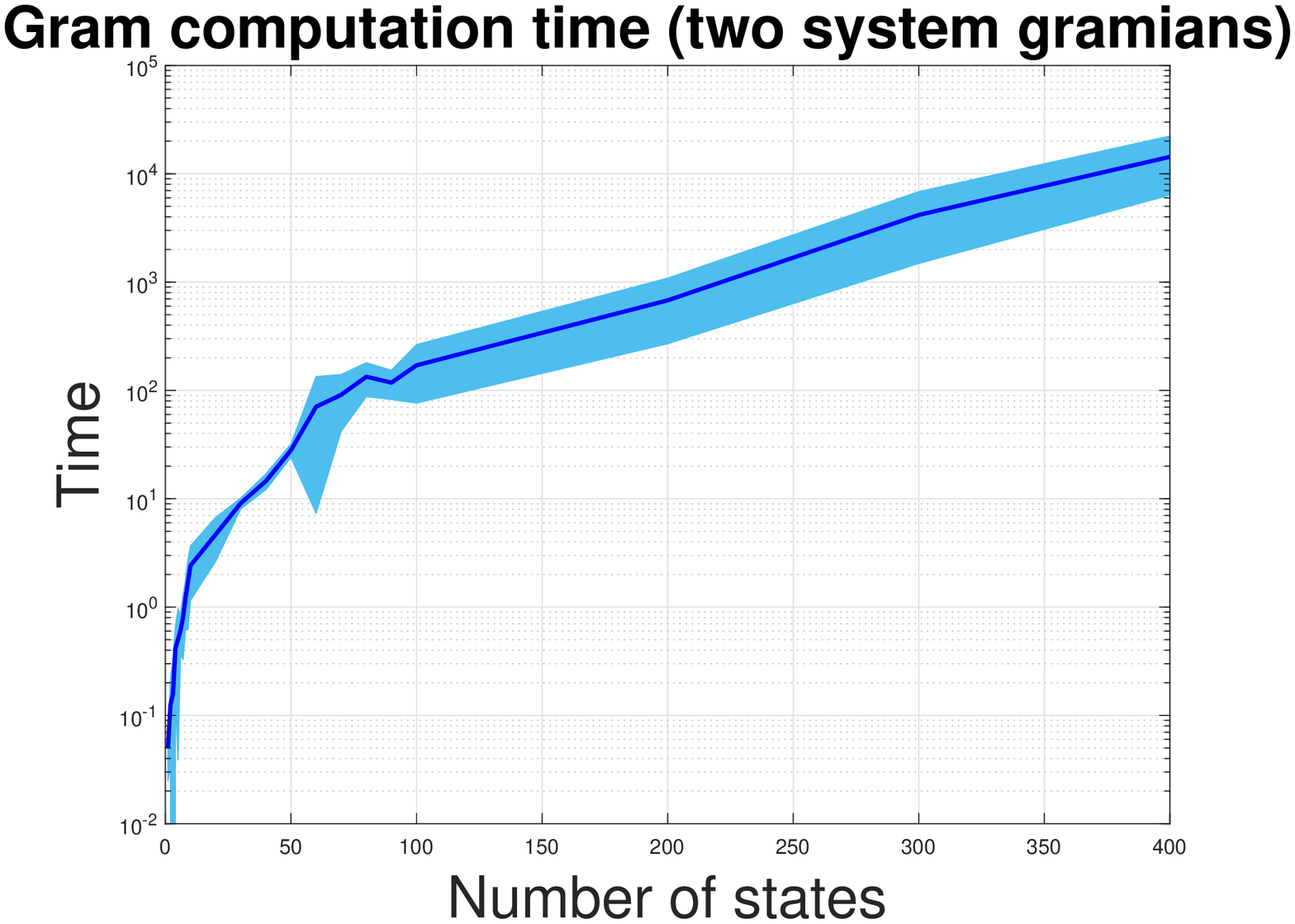}
\end{minipage}
\caption{Average and standard deviation of the computation time of \texttt{h2norm} (left) and \texttt{gram} (right) functions with respect to the number of states. \label{fig:comp}}  
\end{figure}

Figure~\ref{fig:comp} (left) shows that the $\Htwo$ norm computation time is in the order of seconds upto $100$ states, that of minutes upto $300$ states and around $8$ minutes for $500$ states. 

Figure~\ref{fig:comp} (right) shows that the \verb"gram(sys,'co')" computation time is under a minute upto $50$ states, under $10$ minutes upto $200$ states and almost $2$ hours for $400$ states. In the light of these experiments, we see that SSD can handle small to medium size state dimensions. 

For large-scale delay systems with sparse matrices, the \verb"sparse" option of \verb"ssdoptions" can be used to speed up the computation. As a large-scale benchmark, we computed the $\Htwo$ norm of HR$4$ example with $10,000$ states from Section~\ref{sec:bench_h2}. We set the absolute tolerance a little high $10^{-3}$ to make the relative error active and set the relative error to $10^{-4}$ as reported in \cite{Michiels:19} and get similar order of time magnitude, around $6$ seconds, compared to their results $4.8$ seconds obtained by their special large-scale approach.
\begin{verbatim}
>> opt = ssdoptions('Sparse',true, ...
                    'RelTol',1e-4, 'AbsTol',1e-3);
>> tic
>> h2 = h2norm(sys_hr4,opt)
>> time=toc

h2 = 0.4357

time = 5.9651
\end{verbatim}

As a final remark, the function \verb"integral" over $[0,\infty]$ interval may be computationally expensive or less accurate when the integrand has highly oscillatory behavior for a large part of the interval. In our experiments, it shows satisfactory performance.

\section{Concluding Remarks} \label{sec:conclusion}
We took a guided tour of a new MATLAB package, SSD, on several examples. SSD seamlessly integrates with the time and frequency domain visualization capabilities of MATLAB while offering new features on model reduction and $\Htwo$ norm computations built on a frequently-used RTDS delay system representation. The list of SSD functions is given in Table~\ref{tbl:ssd_functions}.

\begin{table}[!h]
\begin{center}
\begin{tabular}{|c|c|}
  \hline
  \hline
  Focus & Function names  \\
  \hline
  Visualizations & \verb"bode", \verb"bodemag", \verb"sigma", \verb"nyquist", \verb"step" \\
  Norm & \verb"h2norm"$^\dag$ \\  
  Model Reduction & \verb"gram"$^\dag$, \verb"balreal", \verb"balred" \\    
  Options & \verb"ssdoptions" \\
  \hline
  \hline
\end{tabular}
\newline
$\dag$ Sparse computation is available.
\caption{SSD functions} \label{tbl:ssd_functions}
\end{center}
\end{table}

Our computational analysis shows that SSD is suitable for small to mid-size state dimensions. We hope that SSD's delay system definition of matrices and delays, easy-to-use interface, and provided functionalities will facilitate the analysis of delay systems and be used as a baseline for advanced techniques.

All the examples in this paper are reproducible by our executable notebook, \verb"introduction.mlx". We share two executable notebooks, \verb"benchmarks_model_reduction.mlx" and \verb"benchmarks_h2norm.mlx" with a collection of benchmark problems on model reduction and $\Htwo$ norm computation.  SSD and three notebooks are available at the GitHub repo
\begin{center}
\texttt{https://github.com/gumussoysuat/ssd$ $}.
\end{center}


\begin{ack}
The author thanks Izmir D. Gumussoy for the fruitful discussions.

\end{ack}


\end{document}